\documentclass[10pt,reqno]{amsart}
  \usepackage{geometry}
\usepackage{amsfonts}
\usepackage{amsmath,amssymb,amsfonts,amsthm}

\usepackage[usenames]{color}

\usepackage[utf8]{inputenc}

\usepackage{url}

\theoremstyle{definition}

\newtheorem{theorem}{Theorem}

\newtheorem{lemma}[theorem]{Lemma}

\newcommand{\C}{\mathbb{C}}
\newcommand{\R}{\mathbb{R}}
\newcommand{\Z}{\mathbb{Z}}

\newcommand{\HH}{\mathbb{H}}

\renewcommand{\Im}{\textrm{Im\,}}

\newcommand{\Sp}{\operatorname{Sp}}
\newcommand{\SU}{\operatorname{SU}}

\newcommand{\UU}{\operatorname{U}}
\newcommand{\Spin}{\operatorname{Spin}}
\newcommand{\Sphere}{\operatorname{S}}

\title{On commutative homogeneous vector bundles attached to nilmanifolds}
\author{Roc\'io D\'iaz Mart\'in and Linda Saal} 
\date{}

\begin{document} 
\maketitle

\textsc{Abstract.} 
The notion of Gelfand pair $(G,K)$ can be generalized by considering homogeneous vector bundles over $G/K$ instead of the homogeneous space $G/K$ and matrix-valued functions instead of scalar-valued functions. This gives the definition of commutative homogeneous vector bundles. Being a Gelfand pair is a necessary condition for being a commutative homogeneous vector bundle. In the case when $G/K$ is a nilmanifold having square-integrable representations, a big family of  commutative homogeneous vector bundles was determined in \cite{rocio2}. In this paper, we complete that classification.

\pagestyle{headings}
\pagenumbering{arabic}

\section{Introduction}

A homogeneous space $G/K$ is called commutative, and the pair $(G, K)$ is called a Gelfand pair,
when $G$ is a locally compact 
group, $K$ is a compact subgroup of $G$, and the convolution algebra $L_0^1(G)$
is commutative. Here $L_0^1
(G)$ denotes the Banach algebra of $L^1$
functions on $G$ satisfying $f(kxk'
) = f(x)$
for $x \in G$ and $k, k' \in K$, where the product  is the usual convolution $(f * h)(g) =
\int_G
f(x)h(x^{-1}
g)dx$ on $G$ (where $dx$ denotes the Haar measure on $G$).

By a nilmanifold we mean a differentiable manifold on which a nilpotent Lie group acts transitively.
By a commutative nilmanifold we mean a commutative space $G/K$ such that $G$ is a Lie group and a closed
nilpotent subgroup $N$ of $G$ acts transitively (cf. e.g. \cite{Wolf2}). In that notation, if $G/K$ is simply connected then 
$N$ acts simply transitively on
$G/K$ and 
$G$ is the semidirect product $K\ltimes N$ (cf. [12,
Theorem 4.2] or \cite[Theorem 13.1.6]{Wolf}). It is shown in \cite{BJR} that if $(K\ltimes N,K)$ is a Gelfand pair then $N$ must be abelian or two-step nilpotent.

Those definitions are associated to the commutative property of the algebra $L^1_0(G)$ of scalar-valued functions. For the vector-valued case we have analogous definitions. 

Again, let $G$ be a Lie group and $K$ a compact subgroup of $G$. It is well known that all the homogeneous vector bundles over the homogeneous space $G/K$ are described by taking finite dimensional representations $(\tau, W_\tau)$ of $K$ (cf. \cite[Section 5.2]{Wal}). Indeed, if we consider the equivalence relation over  $G\times W_\tau$ given by $(gk,w)\sim(g,\tau(k)w)$ for  $g\in G$, $w\in W_\tau$ and $k\in K$, the space of equivalence classes $E_\tau$ has an structure of homogeneous vector bundle over $G/K$. Moreover, each homogeneous vector bundle over $G/K$ is isomorphic to $E_\tau$ for some representation $(\tau,W_\tau)$ of $K$ of finite dimension. 

The space of compactly supported smooth sections of the homogeneous vector bundle $E_\tau$ is naturally identified with the space  of functions $u$  on $\mathcal{D}(G, W_\tau)$ such that $u(xk)=\tau(k)^{-1}u(x)$ for $k\in K$ and $x\in G$. 
It follows from the Schwartz kernel theorem that every linear operator $T$, continuous with 
respect to the standard topologies, mapping $\mathcal{D}$-sections of $E_\tau$ into $\mathcal{D}^\prime$-sections of $E_\tau$ and commuting 
with the action of $G$ on $E_\tau$, can be represented in a unique way as a convolution operator 
\begin{equation}\label{ope}
\big(T(u)\big)(g)=u*F(g):=\int_G F(x^{-1}g)u(x) \ dy \qquad \forall g\in G,
\end{equation}
where  $F$ is a distribution, $F\in \mathcal{D}^\prime(G,\mathrm{End}(W_\tau))$, satisfying  
\begin{equation}\label{F bi tau equiv}
F(k_1xk_2)=\tau(k_2)^{-1}F(x)\tau(k_1)^{-1} \qquad \forall k_1,k_2\in K, \ x\in G.
\end{equation} 
In particular, operators $T = T_F$ as in \eqref{ope} with $F$ a function in $L^1(G,\mathrm{End}(W_\tau))$ satisfying  \eqref{F bi tau equiv}
can be composed with 
each other  and  $$T_{F_2}\circ T_{F_1}=T_{F_1*F_2}.$$ 
We call 
 $L^1_{\tau,\tau}(G,\mathrm{End}(W\tau))$ the convolution algebra of integrable matrix-valued functions with  property \eqref{F bi tau equiv}.

Let $(\tau, W_\tau)$ be an 
irreducible unitary representation of $K$. Let $\widehat{K}$ denote the set of equivalence classes of 
irreducible unitary representations of the group $K$.  The homogeneous vector bundle $E_\tau$ is called commutative, and  the triple $(G,K,\tau)$ is also called commutative, when the algebra $L^1_{\tau,\tau}(G,\mathrm{End}(W_\tau))$ is commutative. In particular,  $(G,K)$ is a Gelfand pair when $(G,K,\tau)$ is a commutative triple with $\tau$ the trivial representation of $K$. Here $L^1_{\tau,\tau}(G,\mathrm{End}(W_\tau))=L^1_0(G)$. It is shown in \cite{Fulvio} that if $G/K$ is connected and if there exists $\tau \in\widehat{K}$ such that $(G, K,\tau)$ is a commutative triple then $(G,K)$ is a Gelfand pair. Therefore, in most cases, being a Gelfand pair is a necessary condition to give rise to commutative triples. When $\tau$ is a character of $K$ and the triple $(G,K,\tau)$ is commutative, these cases are also known as twisted Gelfand pair with respect to the character  $\tau$. Finally, we say that $(G,K)$ is a strong Gelfand pair if $(G,K,\tau)$ is commutative for every $\tau\in\widehat{K}$.

In the case where $L^1_{\tau,\tau}(G,\mathrm{End}(W_\tau))$ is a commutative algebra, the spherical analysis consists of the computation of the continuous characters of such convolution algebra and this gives rise to a kind of simultaneous ``diagonalization'' of all the operators $T_F$ (cf. e.g. \cite{rocio, Koranyi}).

From now on we concentrate on homogeneous vector bundles associated to nilmanifolds. 

The classification of Gelfand pairs $(K\ltimes N,K)$ with $N$ nilpotent was completed by E. Vinberg. 
In the notable article \cite{Vinberg} are exhibited all the Gelfand pairs $(K\ltimes N,K)$ that are irreducible and maximal.  
On the one hand, the irreducibility means that  the center $\mathfrak{z}$ of the Lie algebra $\mathfrak{n}$ of $N$ must be $\mathfrak{z}=[\mathfrak{n},\mathfrak{n}]$ and  $K$ acts irreducibly  on $\mathfrak{n}/\mathfrak{z}$. On the other hand, the maximality implies that the pair $(K\ltimes N, K)$ does not have non-trivial central reductions. This means that $\mathfrak{z}$ does not have $K$-invariant subspaces $\mathfrak{s}$ such that if $\tilde{\mathfrak{s}}$ denotes the orthogonal complement of $\mathfrak{s}$ on $\mathfrak{z}$  and $\tilde{N}$ is the simply connected Lie group with Lie algebra $\tilde{\mathfrak{n}}:=\tilde{\mathfrak{s}}\oplus \left(\mathfrak{n}/\mathfrak{z}\right)\simeq \mathfrak{n}/\mathfrak{s}$ then $(K\ltimes\tilde{N}, K)$ is a Gelfand pair (cf. \cite[Section 13.4A, p. 320]{Wolf}).

Firstly, J. Lauret constructs a family of Gelfand pairs $(K\ltimes N, K)$ considering on $N$ a Riemannian structure (cf. \cite{Lauret}). Here $K\ltimes N$ is the group of isometries of $N$.
When $N$ is particularly the Heisenberg group or the euclidean space $\R^n$, the commutative triples were determined in \cite{Fulvio}. For the corresponding (matrix) spherical analysis see that article and also  \cite{rocio}.
When additionally $N$ has square integrable representations, all the commutative triples that come from these Gelfand pairs were determined  in the recent article \cite{rocio2}. 

This article has the aim of completing the classification of commutative homogeneous vector bundles associated to nilmanifolds, that is, commutative triples of the form $(K\ltimes N,K,\tau)$ with $N$ a nilpotent Lie group. We will analyse which commutative triples come from the Gelfand pairs in \cite{Vinberg} that are not included in the list given in \cite{Lauret}.

\bigskip

\textsc{Acknowledgments.} 
We are immensely grateful to Jorge Lauret who gave impulse to this research.

\section{Preliminaries}

We recall that  since we consider $N$ a two-step nilpotent Lie group, its Lie algebra splits, as a vector space, as $\mathfrak{n}=\mathfrak{z}\oplus V$, where  $V$ is an orthogonal complement of the center $\mathfrak{z}$ and $[V,V]\subset \mathfrak{n}$.
The group $N$ acts naturally on  $\mathfrak{n}$  by  the adjoint action $\mathrm{Ad}$. Also, $N$ acts on $\mathfrak{n}^*$, the real dual space of $\mathfrak{n}$, by the dual or contragredient representation of the adjoint representation
$\mathrm{Ad}^*(n)\lambda:=\lambda\circ \mathrm{Ad}(n^{-1})$  for  $n\in N$ and  $\lambda\in\mathfrak{n}^*$. 
 For each $\lambda\in\mathfrak{n}^*$, let $O(\lambda):=\{\mathrm{Ad}^*(n)\lambda \mid \  n\in N\}$ be its coadjoint orbit. 
From Kirillov's theory there is a correspondence between $\widehat{N}$ and the set of coadjoint orbits. Let  $B_\lambda$ be the skew symmetric bilinear form on $\mathfrak{n}$ given by $$B_\lambda(X,Y):=\lambda([X,Y]) \qquad \forall X,Y\in \mathfrak{n}.$$ 
Let
$\mathfrak{m}\subset\mathfrak{n}$ be a maximal isotropic subalgebra in the sense that $B_\lambda(X,Y)=0$ for all $X,Y\in \mathfrak{m}$ and let $M:=\exp(\mathfrak{m})$. 
Defining on $M$ the character  $\chi_\lambda(\exp(Y)):=e^{i\lambda(Y)}$ for $Y \in\mathfrak{m}$, the irreducible 
representation  $\rho_\lambda\in\widehat{N}$ associated to $O(\lambda)$ is the induced representation $\rho_\lambda:=\mathrm{Ind}_M^N(\chi_\lambda)$. 

Let $X_\lambda\in\mathfrak{z}$ be the representative of $\lambda_{|_\mathfrak{z}}$ (the restriction of $\lambda$ to $\mathfrak{z}$), that is,  $\lambda(Y)=\langle Y,X_\lambda\rangle$ for all $Y\in \mathfrak{z}$. We can split  $\mathfrak{z}=\mathbb{R}X_\lambda\oplus \mathfrak{z}_\lambda$, where $\mathfrak{z}_\lambda:=\mathrm{Ker}(\lambda_{|_{\mathfrak{z}}})$ is the orthogonal complement of $\mathbb{R}X_\lambda$ in $\mathfrak{z}$.

Let $Z$ be the center of $N$. A representation  $(\rho, H_\rho)\in\widehat{N}$ is said square integrable if its matrix coefficients $\langle u,\rho(x)v\rangle$, for $u,v\in H_\rho$,  are square integrable functions on $N$ module $Z$. 

One can see that $\rho_\lambda\in \widehat{N}$ is a square integrable representation if and only if that $B_\lambda$ is non-degenerate on $V$ or, equivalently,  if and only if  the orbits are maximal (cf. \cite{rocio2}).  In this situation, consider the Heisenberg algebra $\mathfrak{n}_\lambda:=\mathbb{R}X_\lambda\oplus V$ with the bracket given by   $$[X,Y]_{\mathfrak{n}_\lambda}:=B_\lambda(X,Y)  X_\lambda.$$ 
Since the character $\chi_\lambda$ is trivial on $\mathfrak{z}_\lambda$, the representation $\rho_\lambda$ acts trivially on $\mathfrak{z}_\lambda$ and defines an irreducible unitary representation of the corresponding Heisenberg group $N_\lambda$. 

Now, let  $P(\lambda)$ be the square root of the determinant of $(B_\lambda)_{|_{V\times V}}$ which is called the Pfaffian.
This function $P$ depends only on $\lambda_{|_{\mathfrak{z}}}$ and so there is a homogeneous polynomial function (which we also denote $P$) on $\mathfrak{z}^*$ such that $P(\lambda)=P(\lambda_{|_{\mathfrak{z}}})$ (cf. \cite[p. 333]{Wolf}).  
According to \cite[Theorem 14.2.10]{Wolf} there is a correspondence between the coadjoint orbits $O(\lambda)$ with $P(\lambda_{|_{\mathfrak{z}}})\neq 0$ and the square
integrable representations.

Apart from that, let $K$ be a compact subgroup of automorphisms of $N$  and let $K_{{\lambda}}$ be the stabilizer of $X_\lambda$ with respect to the action of $K$ on $\mathfrak{n}$. Note that since we always assume $N$ simply connected, we
make no distinction between automorphisms of $N$ and $\mathfrak{n}$. It can be see that $K_\lambda$ is a subgroup of the symplectic group $\mathrm{Sp}(V,(B_\lambda)_{|_{V\times V}})$. Moreover, since $K_\lambda$ is compact, we can assume that it is a subgroup of the unitary group $\mathrm{U}(m)\subset \mathrm{Sp}(V,(B_\lambda)_{|_{V\times V}})$.
Let ${\omega}$ be the the metaplectic representation of $K_\lambda$   associated to the Heisenberg group $N_\lambda$. That is, $$(\omega(k)(p))(z):=p(k^{-1}z)$$ 
for all $p$ the space $\in\mathcal{P}(\mathbb{C}^m)$ of polynomials on $\mathbb{C}^m$, where $2m$ is the dimension on $V$. (For more details see \cite{rocio2}).

For the following theorem see \cite[Theorem 3]{rocio2} and \cite[Theorem 6.1]{Fulvio}. An important fact in the proof is that the classes in $\widehat{N}$ of square integrable representations have full Plancherel measure. 

\begin{theorem}\label{heisenberg}
Let $N$ be a connected and simply connected real two-step nilpotent Lie group which has a square integrable representation. Let $K$ be a compact subgroup of orthogonal automorphisms of $N$ and let $(\tau,W_\tau)\in\widehat{K}$. Then 
$(K\ltimes N,K,\tau)$ is a commutative triple if and only if  $(K_\lambda\ltimes N_\lambda, K_\lambda, \tau_{|_{K_{\lambda}}})$ is a commutative triple for every square integrable representation  $\rho_\lambda\in \widehat{N}$, where $\tau_{|_{K_{\lambda}}}$ denotes the restriction of $\tau$ to $K_{\lambda}$. Also, $(K_\lambda\ltimes N_\lambda, K_\lambda, \tau_{|_{K_{\lambda}}})$ is a commutative triple if and only if ${\omega}\otimes (\tau_{|_{K_{{\lambda}}}})$ is multiplicity free.
\end{theorem}

The previous result will mark the way of our proofs. It allows to make a reduction to Heisenberg groups simplifying our problem from two-step nilpotent Lie groups in general to Heisenberg groups. But it is extremely important the condition on $N$ of having square integrable representations. 

There are few Gelfand pairs $(K\ltimes N,K)$ in the list of E. Vinberg such that the nilpotent group is not of the form given by J. Lauret. Specifically, they correspond to the items 3, 5, 11, 12, 20 and 26 of table 3 of \cite{Vinberg}.

In the items 3 and 26 of \cite{Vinberg} 
the group $N$ does not have square integrable representations. We will skip them from our study. 
In the following paragraphs we will develop the analysis of the triples derived from the remaining Gelfand pairs.
\begin{itemize}
\item[-] Case A: We will develop item 5 of table 3 of \cite{Vinberg}. 
Here we have $\mathfrak{n}=\left[\Lambda^2(\C^{2n})\oplus i\R\right]\oplus \C^{2n}$, where $\Lambda^2(\C^{2n})$ denotes the space of antisymmetric bilinear forms on $\C^{2n}$ over the complex field, and $K=\UU(2n)$.
\item[-] Case B: We will study item 12 of table 3 of \cite{Vinberg}. 
Denoting by $\HH$ the quaternions,  we have $\mathfrak{n}=\left[ H_0(\HH^{n})\oplus Im(\HH)\right]\oplus \HH^{n}$, where $H_0(\HH^{n})$ denotes the space of hermitian $n\times n$ matrices over $\HH$ of trace zero and $Im(\HH)$ the imaginary quaternions, and $K=\Sphere^1\times \Sp(n)$, where $\Sphere^1$ is the one dimensional torus and $\Sp(n)$ the symplectic group. The only difference between item 12 and 11 of table 3 of \cite{Vinberg} is that for item 11 the group $K$ is smaller. 
\item[-] Case C: We will analyze  item 20 of table 3 of \cite{Vinberg}. 
Here we have $\mathfrak{n}= \R^{7}\oplus\R^8$ and $K=\Spin(7)$.
\end{itemize}

\section{Analysis of commutavive and non-commutative triples}\label{section Vinberg}


\subsection*{Case A}

The objects that we will describe in this paragraph can be found in item 5 of table 3 in \cite{Vinberg} as well as in item 5
of table 13.4.1 in the book \cite{Wolf}. 


Consider the two-step nilpotent Lie algebra
$$\mathfrak{n}=\mathfrak{z}\oplus V=\left[\Lambda^2(\C^{2n})\oplus i\R\right]\oplus \C^{2n}$$
where the composition  is given by $$[(z,u),(w,v)]:=((u\wedge v,Im(u^*v)),0) \qquad \forall z, w \in\mathfrak{z}, \ u, v \in V,$$  where we denote by $u^t$ the row vector obtained from the column vector $u$ and by $u^*$ the conjugate row vector. The wedge product $u\wedge v$ can also expressed as $uv^t-vu^t$. 

The unitary group $K=\UU(2n)$ acts on $\mathfrak{z}$ by 
\begin{equation*}
k\cdot(uv^t-vu^t,s):=(k(uv^t-vu^t)k^t,s) \qquad  \forall k\in K \text{ and }  \forall \ uv^t-vu^t\in \Lambda^2(\C^{2n}), \ s\in i\R.
\end{equation*}

\smallskip

We consider in particular a linear functional $\lambda$ in $\mathfrak{n}^*$  such that (when we restrict it to $\mathfrak{z}$)  has as a representative  the element
 $X_\lambda:=(u\wedge v, 0)\in\mathfrak{z}$ with $u\wedge v$ non-degenerate. Let  $\rho_\lambda$ denote the irreducible representation class associated to it. From the formula of the Pfaffian given in  \cite[p. 339--340 ]{Wolf}, $P(\lambda)$ is a positive multiple of the determinant of the matrix $uv^t-vu^t$. According to our election, it is invertible. Therefore $\rho_\lambda$ is square integrable and we are allowed to apply Theorem \ref{heisenberg}. In this case, the subgroup $K_\lambda$ coincides with the symplectic group $\Sp(n)$. 
 
 Let $\tau\in \widehat{K}$ non-trivial. Then there  is a non-trivial $\eta\in\widehat{\Sp(n)}$ appearing in the restriction of $\tau$ to $K_\lambda$.

Now we introduce some notation. Let $\delta_{i,j}$ denote the Kronecker delta. Writing 
$H_i:=\delta_{i,i}-\delta_{n+i,n+i}$, we have that the Cartan subalgebra $\mathfrak{h}$ of $\mathfrak{sp}(n)$ is a complex vector space generated by $\{ H_1,...,H_n \}$.
Let 
$\{L_1,...,L_n\}$
be its dual basis in the dual space $\mathfrak{h}^*$, so $\langle L_i,H\rangle=h_i$ for all $H\in\mathfrak{h}$.  From the theorem of the highest weight,  every  irreducible representation of $\mathfrak{sp}(n)$ is in correspondence with a non-negative integer linear combination of the fundamental weights. Hence $\eta \in \widehat{\mathrm{Sp}(n)}$ can be parametrized in terms of the weights $\{L_i\}$ as $(\eta_1,...,\eta_n)$  where  $\eta_i\in\mathbb{Z}_{\geq 0}$ $\forall i$ and $\eta_1\geq \eta_2\geq ... \geq \eta_n$. (For a reference see, for example,  \cite{Fulton y Harris, Knapp}.) We will denote the representation $\eta\in\widehat{\mathrm{Sp}(n)}$ by $\eta_{(\eta_1,...,\eta_n)}$ to emphasize that the representation $\eta$ is in correspondence with the tuple $(\eta_1,...,\eta_n)$ that we call partition. 
 
The metaplectic representation  $(\omega, \mathcal{P}(\C^{2n}))$ of $\Sp(n)$ decomposes (using the notation given in the previous sections) as 
\begin{equation*}
\omega=\bigoplus_{j\in\Z_{\geq 0}} \eta_{(j)},
\end{equation*}
where each
$\eta_{(j)}$ corresponds to the partition on length one $(j)$, for some non-negative integer $j$.

The following fact is proved in \cite[Corollary 2]{rocio2} and will be useful to deduce our main results.
 
\begin{lemma}
\label{coro 1 sp} Let $\eta$ be an irreducible representation of $\mathrm{Sp}(n)$ then
$\eta$ appears in the decomposition into irreducible factors of $\eta\otimes\eta_{(2)}$.
\end{lemma}

From the above lemma, the representation $\eta$ 
 appears in the decomposition of the factors $\eta\otimes\eta_{(0)}$ and $\eta\otimes\eta_{(2)}$. Then, it
appears with multiplicity in the decomposition into irreducible factors of $\omega\otimes\tau_{|_{K_{\lambda}}}$. 
Therefore, by Theorem \ref{heisenberg}, we have the following result.

\begin{theorem}\label{th 1}
The triple $(K\ltimes N, K, \tau)$, where $N$ is the simply connected nilpotent Lie group with Lie algebra $\left[\Lambda^2(\C^{2n})\oplus i\R\right]\oplus \C^{2n}$ and  $K=\UU(2n)$, is commutative if and only if $\tau$ is the trivial representation.
\end{theorem}

\subsection*{Case B}

The objects that we will describe in this paragraph can be found in item 12 of table 3 in \cite{Vinberg} as well as in item 9
of table 13.4.1 in the book \cite{Wolf}.



Consider the two-step nilpotent Lie algebra
$$\mathfrak{n}=\mathfrak{z}\oplus V:=\left[ H_0(\HH^{n})\oplus \Im(\HH)\right]\oplus \HH^{n}$$ with Lie bracket given by $$[(z,u),(w,v)]:=(((uiv^*-viu^*)_0,u^*v-v^*u)),0)$$ 
 for all $z, w \in\mathfrak{z}= H_0(\HH^{n})\oplus \Im(\HH)$ and  $ u, v \in V=\HH^{n}$, where, in general, if $A$ is a matrix, $(A)_0:=A-\frac{1}{2}tr(A)I$. Let $N$ be the simply connected nilpotent Lie group with Lie algebra $\mathfrak{n}$. 

Every point in group of automorphisms $(e^{i\theta},k)\in K= \Sphere^1\times \Sp(n)$  acts on $\mathfrak{n}$ by 
\begin{equation*}
(e^{i\theta},k)\cdot((A,q),v):=((kAk^*,e^{i\theta} q e^{-i\theta}),kv e^{-i\theta}),\end{equation*}
for all $A\in H_0(\HH^{n})$, $q\in Im(\HH)$ and  $v\in\HH^{n}$.
\smallskip

Let $\tau=\eta\otimes\chi_r$ where $\eta\in\widehat{\Sp(n)}$ and $\chi_r$ a character of $\Sphere^1$.  

First we consider the case where the factor in $\widehat{\Sp(n)}$ is trivial, that is $\tau=\chi_r$. Since the convolution algebra $L_{\tau,\tau}^1(K\ltimes N,\mathrm{End}(W_\tau))$ is  naturally identified with the space of $\mathrm{End}(W_\tau)$-valued integrable functions on $N$ such that $F(k\cdot x) =\tau(k)F(x)\tau(k)^{-1}$ for all $k\in K$ and $x\in N$ (cf. e.g. \cite{rocio,Fulvio}). In this situation, it totally coincides with the convolution algebra of $K$-invariant integrable scalars functions on $N$, which is commutative since we have a Gelfand pair. Therefore we have commutative triples. 

Now let $\eta\in\widehat{\Sp(n)}$ non-trivial. 

We consider the linear functional $\lambda$ with representative 
$X_\lambda:=(0, q)\in\mathfrak{z}$. Let  $\rho_\lambda$ denote the irreducible representation class associated to it. 
From the formula of the Pfaffian given \cite[p. 340]{Wolf}, $P(\lambda)$ is a positive multiple of $|q|^{2n}$. Then for all imaginary non-zero quaternions $q$, the associated representation  $\rho_\lambda$ is square integrable. The group $K_\lambda$ is easily calculated: If the quaternion $q$ belongs to $i\R$, $K_\lambda$ coincides with $K$ (since $q$ commutes with every complex); for the other cases $K_\lambda$ is $\Sp(n)$. We consider for example $q=j\in Im(\HH)$ in order to fix $K_\lambda=\Sp(n)$. Note that  $\tau_{|_{{K_{\lambda}}}}=\eta$.
 
  The metaplectic representation  $(\omega, \mathcal{P}(\C^{2n}))$ of $K_\lambda=\Sp(n)$ decomposes as 
 \begin{equation*}
 \omega=\bigoplus_{j\in\Z_{\geq 0}} \eta_{(j)}.
 \end{equation*}
 Therefore  $\omega\otimes\eta$ is not multiplicity free: $\eta$ appears in the  the factors $\eta\otimes\eta_{(0)}$ and in $\eta\otimes\eta_{(2)}$ by Lemma \ref{coro 1 sp}.
Then, by Theorem \ref{heisenberg}, in these cases the triples are not  commutative.

\begin{theorem}\label{th 2}
The triple $(K\ltimes N, K, \tau)$, where $N$ is the simply connected nilpotent Lie group with Lie algebra $\left[ H_0(\HH^{n})\oplus Im(\HH)\right]\oplus \HH^{n}$ and  $K=\Sphere^1\times \Sp(n)$, is commutative if and only if $\tau\in\widehat{\Sphere^1}$.
\end{theorem}

To conclude this case we want to note that the analysis of the item 11 of table 3 of \cite{Vinberg} is the same (or simpler) since the subgroup $K$ is only $\Sp(n)$.

\subsection*{Case C}

The objects that we will describe in this paragraph can be found in item 20 of table 3 in \cite{Vinberg} as well as in item 2
of table 13.4.1 in the book \cite{Wolf}.

In this case we study an H-type group (\cite{Kaplan}). Consider the two-step nilpotent Lie algebra $\mathfrak{n}=\mathfrak{z}\oplus V$ where  $V$ is $\R^8$ or the octonions $\mathbb{O}$ and the center $\mathfrak{z}$ is $\R^7$  or the imaginary octonions $Im(\mathbb{O})$ with Lie bracket on $V$ characterized by:
${\langle[u, v], z \rangle}_\mathfrak{z}:={\langle J(z)u, v\rangle}_V$,
where $(J,V)$
is a real representation of $\mathfrak{z}$ given by the product on the octonions $J(z)v:=zv$ for  $z\in Im(\mathbb{O})$ and $v\in \mathbb{O}$. 

Let $K$ be the maximal connected group of orthogonal automorphisms of $N$. Precisely, $K=\Spin(7)$.
(We mention that in general, for an H-type group, it group of automorphisms was determined by L. Saal in  \cite{la linda}.)

Let $\lambda$ be a linear functional  with representative 
$X\in\mathfrak{z}$. From  \cite[p. 340]{Wolf}, the Pfaffian $P(\lambda)$ is not null almost everywhere. In this case, $K_\lambda$ is isomorphic to the spin group $\Spin(6)$. Therefore we will work with the Lie algebra $\mathfrak{so}(6)$ .

Let $\tau$ be a non-trivial irreducible unitary representation of $\Spin(7)$. When we restrict $\tau$ to $\Spin(6)$ it decomposes as a sum of irreducible representations  and we can pick one non-trivial factor. We will denote it as  $\tilde{\eta}\in\widehat{\Spin(6)}$ and we will identify it with its derived representation and also we will view it as a representation of $\mathfrak{so}(6)$. 

Apart from that, it is easy to derive the decomposition of the metaplectic representation since $\Spin(6)$ is isomorphic to $\SU(4)$ and it is well known that its action on $\mathcal{P}(\mathbb{C}^4)$ decomposes into $\oplus_{j\in\Z_{\geq 0}}\mathcal{P}_j(\mathbb{C}^4)$. It can be identified with 
\begin{equation}\label{omega}
 \omega=\bigoplus_{j\in\Z_{\geq 0}} \tilde{\eta}_{(j)}.
 \end{equation}
where $\tilde{\eta}_{(j)}$ corresponds to the partition of length one $(j)$ of $\mathfrak{so}(6)$. Here we decided to use an analogous notation  to that of the representations of $\mathfrak{sp}(n)$. 
As a consequence of the highest weight theorem, 
we can associate each irreducible representation $\sigma$ of $\mathfrak{so}(2m)$ with an $n$-tuple of integers  $(\sigma_1,...,\sigma_m)$ satisfying the condition $\sigma_1\geq  \sigma_2\geq...\geq\sigma_{m-1}\geq |\sigma_m|$ (cf. \cite{Fulton y Harris, Okada}). The following fact can be found in \cite[Theorem 3.2 and Remark 3.5]{Okada}.

\begin{lemma}\label{Okada} 
Let $\tilde{\eta}$ be an arbitrary irreducible representation of $\mathfrak{so}(2m)$ associated to a sequence of $m$ integers $\tilde{\eta}_1\geq \tilde{\eta}_2\geq...\geq \tilde{\eta}_{m-1}\geq |\tilde{\eta}_m|$ and let $s$ be a non-negative integer.   Then the multiplicity of an irreducible representation  $\sigma$ (associated to the sequence of integers $\sigma_1\geq \sigma_2\geq...\geq \sigma_{m-1}\geq |\sigma_m|$) in the tensor product $\tilde{\eta}\otimes\tilde{\eta}_{(s)}$ is equal to the number of integer sequences $\varsigma$ satisfying:
\begin{itemize}
    \item[$(i)$] $\varsigma_{1} \geq \varsigma_{2} \geq ... \geq \varsigma_{m-1} \geq\left|\varsigma_{m}\right|$
    \item[$(ii)$] $\tilde{\eta}_{1}\geq\varsigma_{1} \geq \tilde{\eta}_{2}\geq\varsigma_{2} \geq ... \geq\varsigma_{m-1}\geq \tilde{\eta}_{m}\geq\varsigma_{m}$ and $\sigma_{1}\geq\varsigma_{1} \geq \sigma_{2}\geq\varsigma_{2} \geq ... \geq\varsigma_{m-1}\geq \sigma_{m}\geq\varsigma_{m}$
    \item[$(iii)$] $\sum_{i=1}^{m}\left(\tilde{\eta}_{i}-\varsigma_{i}\right)+\sum_{i=1}^{m}\left(\sigma_{i}-\varsigma_{i}\right)=s$
    \item[$(iv)$] $\varsigma_m\in\{\tilde{\eta}_m,\sigma_m\}$
\end{itemize}
\end{lemma}

Now will apply this lemma, with $m=3$ a fixed non-trivial irreducible representation $\tilde{\eta}$ in the decomposition of $\tau_{|_{K_{\lambda}}}$,  to deduce whether  the tensor product $\omega\otimes\tau_{|_{K_{\lambda}}}$ is  multiplicity free or not. 
The representation $\tilde{\eta}$ corresponds to $(\tilde{\eta}_1,\tilde{\eta}_2,\tilde{\eta}_3)$ with $\tilde{\eta}_i\in\mathbb{Z}$ for $i=1,2,3$ and $\tilde{\eta}_1\geq \tilde{\eta}_2\geq |\tilde{\eta}_3|$ and $\tilde{\eta_1}>0$ since it is non-trivial. 
\begin{itemize}
    \item If  $\tilde{\eta}_1>\tilde{\eta}_3$ (in particular if $\tilde{\eta}_3\leq 0$), then $\tilde{\eta}$ appears in $\tilde{\eta}\otimes \tilde{\eta}_{s}$ for $s=2(\tilde{\eta}_1-\tilde{\eta}_3)>0$ since  $\varsigma=(\tilde{\eta}_2,\tilde{\eta}_3,\tilde{\eta}_3)$ satisfies the conditions listed in Lemma \ref{Okada}. By (\ref{omega}) we have that $\tilde{\eta}$ appears at least twice in $\omega\otimes\tau_{|_{K_{\lambda}}}$.
    \item  If  $\tilde{\eta}_1=\tilde{\eta}_2=\tilde{\eta}_3=r$ for $r\in\Z_{>0}$, then if $\sigma=(\sigma_1,\sigma_2,\sigma_3)$ appears in $\tilde{\eta}\otimes \tilde{\eta}_{s}$ for $s\in\Z_{\geq 0}$, then by Lemma \ref{Okada}, $\sigma_1=\sigma_2=r$. Also, with the notation in the Lemma, it must happen that $\varsigma$ is $(r,r,r)$ or $(r,r,\sigma_3)$. For the first case $\sigma_3-r=s$, then since $r\geq\sigma_3$  
    the only possibility for $s$ is $s=0$.
    For the second, $r-\sigma_3=s$. Therefore $s$ is completely determined by $\sigma$. Thus $\sigma=(r,r,\sigma_3)$ with $r\geq \sigma_3$ appears once in  $\tilde{\eta}\otimes \oplus_{j\in\Z_{\geq 0}}\tilde{\eta}_{(j)}$. Consequently, if $\tilde{\eta}=(r,r,r)$ and $\sigma=(r,r,\sigma_3)$ with $r> \sigma_3$ both appears in $\tau_{|_{K_{\lambda}}}$, the tensor product $\omega\otimes\tau_{|_{K_{\lambda}}}$ is not  multiplicity free. Note that in this situation $\sigma$ satisfies the conditions $\tilde{\eta}$ for the first item. Therefore we have the following conclusion.   
\end{itemize}

\begin{theorem}\label{th 3}
The triple   $(K\ltimes N, K, \tau)$, where $N$ is the simply connected nilpotent Lie group with Lie algebra $Im(\mathbb{O})\oplus\mathbb{O}$ and  $K=\Spin(7)$, is commutative if and only if $\tau_{|_{\Spin(6)}}$ is associated to a partition of the form $(r,r,r)$ for  $r\in\Z_{\geq 0}$.
\end{theorem}

\section{Conclusions}

At this point we sum up the following.  Theorem \ref{th 1} does not provide non-trivial commutative triples. Theorem \ref{th 2} gives rise to a commutative triple only if $\tau$ is a character of $\Sphere^1$. These cases are 
twisted Gelfand pairs. 

Finally, Theorem \ref{th 3} gives a family of commutative triples. This family is similar to the case of the H-type group listed by J. Lauret in \cite{Lauret}. There the nilpotent Lie algebra $\mathfrak{n}$ is $Im(\mathbb{H})\oplus\mathbb{H}^n$ 
where $Im(\mathbb{H})$ is its center $\mathfrak{z}$ and the Lie bracket on $\mathbb{H}^n$ is given, analogously as in case C, by the real representation of 
$\mathfrak{z}$, 
$J(z)(v):=(zv_1,...,zv_n)$ for  $v=(v_1,...,v_n)\in\mathbb{H}$. Its (maximal connected) group of automorphisms is $K=\mathrm{SU}(2)\times \mathrm{Sp}(n)$. By \cite[Proposition 1]{rocio2}, in that case we obtain commutative triples if and only if  $\tau\in\widehat{\mathrm{SU}(2)}$ or
  $\tau\in \widehat{\mathrm{Sp}(n)}$
corresponding to a partition of the form $(a,a,...,a)$ of length at most $n$ for a non-negative integer $a$.

In conclusion, denoting by $\mathrm{H}_n$ the Heisenberg group $\mathrm{H}_n=\C_n\times \R$, if we exclude from the analysis the Gelfand pairs of the form $(\mathrm{H}_n\ltimes K,K)$ where $K$ is a proper subgroup of $\mathrm{U}(n)$, all the commutative triples of the form $(K\ltimes N, K,\tau)$, where $N$ is a simply connected nilpotent Lie group having square integrable representations and $K$ is a compact subgroup of $G$, are the following:

\begin{enumerate}

\item $([\Sphere^1\times \Sp(n)]\ltimes N, \Sphere^1\times \Sp(n), \tau)$, where $N$ is the simply connected nilpotent Lie group with Lie algebra $\left[ H_0(\HH^{n})\oplus Im(\HH)\right]\oplus \HH^{n}$, is commutative for $\tau\in\widehat{\Sphere^1}$. It corresponds to the preceding case C. It is a twisted Gelfand pair.

\item $(\mathrm{SU}(n)\times \mathbb{S}^1, N(\mathfrak{su}(n),\mathbb{C}^n),\tau)$, for all $\tau\in\widehat{\mathbb{S}^1}$, where $n\geq 3$. It is a twisted Gelfand pair.

\item $(\mathrm{SU}(n)\times \mathbb{S}^1, N(\mathfrak{u}(n),\mathbb{C}^n),\tau)$, for all $\tau\in\widehat{\mathbb{S}^1}$, where $n\geq 3$. It is a twisted Gelfand pair.

\item $(\Spin(7)\ltimes N, \Spin(7), \tau)$, where $N$ is the simply connected nilpotent Lie group with Lie algebra $Im(\mathbb{O})\oplus\mathbb{O}$, is commutative when $\tau_{|_{\Spin(6)}}$ is associated to a constant partition $(r,r,r)$ for  $r\in\Z_{\geq 0}$. $N$ is of Heisenberg type. It corresponds to the preceding case C.

\item $\left(\mathrm{SU}(2)\times \mathrm{Sp}(n), N(\mathfrak{su}(2),(\mathbb{C}^2)^n),\tau\right)$, for all $\tau\in\widehat{\mathrm{SU}(2)}$ and for all $\tau\in\widehat{\mathrm{Sp}(n)}$ associated to a constant partition of length at most $n$, where $n\geq 1$. $N(\mathfrak{su}(2),(\mathbb{C}^2)^n)$ is of Heisenberg type.

\item $(\mathrm{SU}(2)\times \mathrm{U}(k)\times \mathrm{Sp}(n), N(\mathfrak{u}(n),(\mathbb{C}^2)^k\oplus(\mathbb{C}^2)^n),\tau)$,   for all $\tau\in \widehat{\mathrm{U}(k)}$, where $k\geq 1, n\geq 0$.

\item $(G\times U, N(\mathfrak{g}, V), \tau)$, where $
 \mathfrak{g}=\mathfrak{su}(m_1)\oplus...\oplus\mathfrak{su}(m_\beta)\oplus\mathfrak{su}(2)\oplus...\oplus\mathfrak{su}(2)\oplus\mathfrak{c}$ 
 with  $\mathfrak{c}$  an abelian component,  
$G=\mathrm{SU}(m_1)\times...\times \mathrm{SU}(m_\beta)\times \mathrm{SU}(2)\times...\times \mathrm{SU}(2)$, $
V= \mathbb{C}^{m_1}\oplus...\oplus \mathbb{C}^{m_\beta}\oplus\mathbb{C}^{2k_1+2n_1}\oplus...\oplus\mathbb{C}^{2k_\alpha+2n_\alpha}$  and $U=\mathbb{S}^1\times...\times \mathbb{S}^1\times \mathrm{U}(k_1)\times \mathrm{Sp}(n_1)\times...\times \mathrm{U}(k_\alpha)\times \mathrm{Sp}(n_\alpha)$,
for all $\tau\in\widehat{\mathbb{S}^1}\otimes...\otimes \widehat{\mathbb{S}^1}\otimes  \widehat{\mathrm{U}(k_1)}\otimes...\otimes \widehat{\mathrm{U}(k_\alpha)}$, where $m_j\geq 3$ for all $1\leq j\leq\beta$, $k_i\geq 1$, $n_i\geq 0$ for all $1\leq i\leq\alpha$.

\item $(\mathrm{H}_n\ltimes \mathrm{U}(n), \mathrm{U}(n),\tau)$ for all $\tau\in\widehat{\mathrm{U}(n)}$ (proved by Yakimova in \cite{Yakimova}). This is the unique strong Gelfand pair of this form.
\end{enumerate}

Where in \'items 2, 3, 5, 6 and 7 the nilpotent Lie group $N=N(\mathfrak{g},V)$ is 
endowed with a left-invariant Riemannian metric determined by an inner product $\langle\cdot,\cdot\rangle$ on its Lie algebra $\mathfrak{n}$ described as follows. 
Let $(\pi, V)$ be  a faithful real representation of a compact Lie algebra $\mathfrak{g}$. We consider  inner products $\langle\cdot,\cdot\rangle_\mathfrak{g}$ on $\mathfrak{g}$ and $\langle\cdot,\cdot\rangle_V$  on $V$ such that $\langle\cdot,\cdot\rangle_\mathfrak{g}$  is   $\mathrm{ad}(\mathfrak{g})$-invariant and $\langle\cdot,\cdot\rangle_V$ is $\pi(\mathfrak{g})$-invariant.
Let $\mathfrak{n}:=\mathfrak{g}\oplus V$ be
the two-step nilpotent Lie algebra with center $\mathfrak{g}$ and Lie bracket defined on $V$ by
${\langle[u, v], X \rangle}_\mathfrak{g}:={\langle\pi(X)u, v\rangle}_V$ for all $u, v \in V$, $X \in \mathfrak{g}$. These inner products define an inner product $\langle\cdot,\cdot\rangle$ on $\mathfrak{n}$ satisfying $\langle\cdot,\cdot\rangle_{|_{\mathfrak{g}\times \mathfrak{g}}}= \langle\cdot,\cdot\rangle_\mathfrak{g}$,
$\langle\cdot,\cdot\rangle_{|_{V\times V}}= \langle\cdot,\cdot\rangle_V$
 and $\langle\mathfrak{g},V\rangle=0$.  Let $N(\mathfrak{g}, V)$ be the two-step connected simply connected Lie group with Lie algebra $\mathfrak{n}$. Specifically, in item 2, $\mathbb{C}^n$ denotes the standard representation of $\mathfrak{su}(n)$, in item 3, $\mathbb{C}^n$ denotes the standard representation of $\mathfrak{u}(n)$ regarded as a real representation, in item 5, $\mathbb{C}^2$ denotes the standard representation of $\mathfrak{su}(2)$ regarded as a real representation, in item 6, $(\mathbb{C}^2)^k\oplus(\mathbb{C}^2)^n$ is an orthogonal sum, the center of $\mathfrak{u}(2)$ acts non-trivially only on $(\mathbb{C}^2)^k$, $(\mathbb{C}^2)^n$ denotes the representation of $\mathfrak{su}(2)$ stated in the first item,  $\mathfrak{u}(2)$ acts component-wise on $(\mathbb{C}^2)^k$ in the standard way regarded as a real representation and $\tau\in \widehat{\mathrm{U}(k)}$ and in item 7, $\mathfrak{g}$ is acting on $V$ as follows: For each $1\leq i\leq \beta+\alpha$, $\mathfrak{c}$ has a unique subspace $\mathfrak{c}_i$   acting non-trivially on only the $i$-th component of $V$ and the dimension of $\mathfrak{c}_i $ is $1$.  For $1\leq i\leq \beta$,
 $\mathfrak{su}(m_i)\oplus \mathfrak{c}_i$ (which is isomorphic to $\mathfrak{u}(m_i)$)  acts non-trivially only on  $\mathbb{C}^{m_i}$ by standard. For $\beta+1<i\leq \beta+\alpha $, $\mathfrak{su}(2)\oplus\mathfrak{c}_i$  acts non-trivially only on $\mathbb{C}^{2k_i+2n_i}$ as in the above case. For more details see \cite{rocio2, Lauret, Lauret nil}.

\end{document}